\documentclass{amsart}
\usepackage{epsfig}

\newtheorem{theorem}{Theorem}[section]
\newtheorem{lemma}[theorem]{Lemma}
\newtheorem{proposition}[theorem]{Proposition}
\newtheorem{corollary}[theorem]{Corollary}

\theoremstyle{definition}
\newtheorem{definition}[theorem]{Definition}

\theoremstyle{remark}

\begin{document}

\title[$2$-adic valuation]{The $2$-adic valuation of a 
sequence arising from a rational integral}

\author{Tewodros Amdeberhan}
\address{Department of Mathematics,
Tulane University, New Orleans, LA 70118}
\email{tamdeberhan@math.tulane.edu}

\author{Dante Manna}
\address{Department of Mathematics and Statistics, Dalhousie University, 
Halifax, Nova Scotia, canada, B3H 3J5}
\email{dmanna@mathstat.dal.ca}

\author{Victor H. Moll}
\address{Department of Mathematics,
Tulane University, New Orleans, LA 70118}
\email{vhm@math.tulane.edu}


\date{\today}


\begin{abstract}
We analyze properties of the $2$-adic valuations of an integer
sequence that originates from an explicit 
evaluation of a quartic integral. We also give a combinatorial 
interpretation of the valuations of this sequence. Connections with 
the orbits arising from the Collatz problem are discussed. 
\end{abstract}

\maketitle

\newcommand{\realpart}{\mathop{\rm Re}\nolimits}
\newcommand{\imagpart}{\mathop{\rm Im}\nolimits}

\numberwithin{equation}{section}

\section{Introduction} \label{intro} 
\setcounter{equation}{0}

The sequence
\begin{equation}
A_{l,m}  =  \frac{l! \, m!}{2^{m-l}} 
\sum_{k=l}^{m} 2^{k} \binom{2m-2k}{m-k} \binom{m+k}{m} 
\binom{k}{l} 
\label{intseq}
\end{equation}
\noindent
for $m \in \mathbb{N}$ and $0 \leq l \leq m$ appears in the  
evaluation of the definite integral
\begin{equation}
N_{0,4}(a;m) = \int_{0}^{\infty} \frac{dx}{(x^{4}+4ax^{2} + 1)^{m+1}}. 
\label{nzero4}
\end{equation}
\noindent
Explicitly,
\begin{equation}
N_{0,4}(a;m) = \frac{\pi}{\sqrt{2} \, m! \, (4(2a+1))^{m+1/2}} 
\sum_{l=0}^{m} A_{l,m} \frac{a^{l}}{l!}. 
\end{equation}

The evaluation of $A_{l,m}$ using (\ref{intseq}) is efficient if $l$ is 
close to $m$. For instance,
\begin{equation}
A_{m,m} = 2^{m} (2m)! \text{ and } A_{m-1,m} = 2^{m-1}(2m-1)! (2m+1). 
\end{equation}

In \cite{bomosha} it is shown that $A_{l,m}$ is 
always an integer. An efficient method for the evaluation of these
sequences when $l$ is small is presented there. For example,
\begin{equation}
A_{0,m} = \prod_{k=1}^{m} (4k-1) \text{ and } 
A_{1,m} = (2m+1) \prod_{k=1}^{m}(4k-1) - 
\prod_{k=1}^{m} (4k+1). 
\label{azero}
\end{equation}

The results described in this paper started as empirical observations on 
the behavior of $\nu_{2}(A_{l,m})$, the $2$-adic valuation of $A_{l,m}$. Recall that $\nu_{2}(x)$
is the highest power of $2$ that divides $x$.

The $2$-adic valuation of $A_{0,m}$ follows directly from 
(\ref{azero}). Clearly $A_{0,m}$ is odd, so $\nu_{2}(A_{0,m}) = 0$. The 
$2$-adic valuation of $A_{1,m}$ is given by 
\begin{equation}
\nu_{2}(A_{1,m}) = \nu_{2}(m(m+1)) + 1. 
\label{2value1}
\end{equation}
\noindent
This is the main result of \cite{bomosha}. 

The first goal of this paper is to present the following generalization of 
(\ref{2value1}). 

\begin{theorem}
\label{2adicall}
The $2$-adic valuation of $A_{l,m}$ satisfies 
\begin{equation}
\nu_{2}(A_{l,m}) = \nu_{2}((m+1-l)_{2l}) + l, 
\label{2valuel}
\end{equation}
\noindent
where $(a)_{k} = a(a+1) \cdots (a+k-1)$ is the Pochhammer symbol.
\end{theorem}

As a consequence of this theorem we prove some interesting combinatorial 
properties of the sequence $A_{l,m}$. Henceforth, we assume that
the index $l \in \mathbb{N}$ is fixed and $m \geq l$.

Figure 1 shows the 
graph of $\nu_{2}(A_{60,m})$ for $60 \leq m \leq 450$. The horizontal axis 
is the translate $m' = m - 59$, so the indexing starts at $1$.  \\

{{
\begin{figure}[ht]
\begin{center}
\centerline{\psfig{file=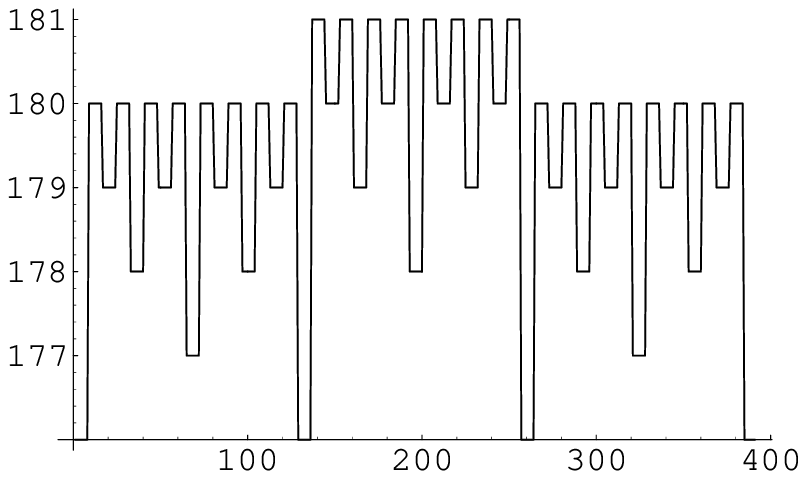,width=35em,angle=0}}
\caption{The $2$-adic valuation of $A_{60,m}$ for $1 \leq m' \leq 400$}
\label{figure1}
\end{center}
\end{figure}
}}

\medskip

The figure suggests that the values of $\{ \nu_{2}(A_{60,m}): m \geq 60 \}$ 
have a {\em block structure} meaning that they are 
composed of consecutive blocks,
all of the same length. Indeed, this sequence begins with
\begin{eqnarray}
& & \{ 
176, \, 176, \, 176, \, 176, \, 176, \, 176, \, 176, \, 176 , \,
 180, \, 180, \, 180, \, 180, \, 180, \, 180, \, 180, \, 180, 
\nonumber \\
& & \, 179, \, 179, \, 179, \, 179, \, 179, \, 179, \, 179, \, 179 , \,
 180, \, 180, \, 180, \, 180, \, 180, \, 180, \, 180, \, 180, \, \ldots 
\nonumber 
\},
\nonumber
\end{eqnarray}
\noindent
which is formed by blocks of length $8$. 

\medskip

This motivates the next definition.

\begin{definition}
Let $s \in \mathbb{N}, \, s \geq 2$. We say that a sequence 
$\{ a_{j}: \, j \in \mathbb{N} \}$ is {\em simple of length} $s$ ( or 
{\em $s$-simple}) if $s$ is the largest integer such that, for 
each $t \in \{ 0, \, 1, \, 2, \, \cdots \}$, we have
\begin{equation}
a_{st + 1} = a_{st + 2} = \cdots = a_{s(t+1)}.
\end{equation}
\noindent
The sequence $\{ a_{j}: \, j \in \mathbb{N} \}$ is said to have a 
{\em block structure} if 
it is $s$-simple for some $s \geq 2$. 
\end{definition}

Using theorem \ref{2adicall} we evaluate
\begin{equation}
\nu_{2}(A_{l,m+1}) - \nu_{2}(A_{l,m})  =  \nu_{2}(m+l+1) - \nu_{2}(m-l+1). 
\end{equation}
\noindent
We then use this fact to establish in Theorem \ref{period-thm} that the 
sequence of integers 
$\{ \nu_{2}(A_{l,m}): m \geq l \}$ is $2^{1 + \nu_{2}(l)}$-simple. 

The combinatorial properties of the sequence $\nu_{2}(A_{l,m})$ are 
described as an algorithm:  \\

\noindent
{\bf The maps $F$ and $T$}. Consider the operators defined on sequences by: 
\begin{eqnarray}
F ( \{ a_{1}, \, a_{2}, \, a_{3}, \, \cdots \} )  & := & 
\{ a_{1}, \, a_{1}, \, a_{2}, \, a_{3}, \, \cdots \}, \label{Fdef}
\end{eqnarray}
\noindent
and
\begin{eqnarray}
T ( \{ a_{1}, \, a_{2}, \, a_{3}, \, \cdots \} ) & := & 
\{ a_{1}, \, a_{3}, \, a_{5}, \, a_{7}, \, \cdots \}. \label{Tdef}
\end{eqnarray}
\noindent
Now introduce the sequence $c$ as
\begin{eqnarray}
c & := & \{ \nu_{2}(m): \quad m \geq 1 \}
= \{ 0, \, 1, \, 0, \, 2, \, 0, \, 1, \, 0, \, 3, \, 0, \cdots \}.
\label{cdef}
\end{eqnarray}

\medskip

\noindent
{\bf The algorithm}: \\

\noindent
1) Start with the sequence
$X(l)  :=  \left\{ \nu_{2}( A_{l}(l+m-1)): \quad m \geq 1 \, \right\}. $

\medskip

\noindent
2) Find $n \in \mathbb{N}$ so that the sequence $X(l)$ 
is $2^{n}$-simple. Define 
$Y(l)  :=  T^{n} \left( X(l) \right)$. At the initial stage, 
Theorem \ref{period-thm}
ensures that $n= 1 + \nu_{2}(l)$. 

\medskip

\noindent
3) Introduce the shift 
$Z(l) :=   Y(l) - c$.  

\medskip

\noindent
4)  Define $W(l):= F(Z(l))$. 

\medskip

\noindent
If $W$ is a constant sequence, then STOP; otherwise go to step $2)$ 
with $W$ instead of $X$.  Define $X_k(l)$ as the new 
sequence at the end of the $(k-1)$th cycle of this process, with
$X_1(l)=X(l)$.

\medskip

The next theorem justifies that the steps 
described above make sense and comprise an algorithm. In 
Section \ref{verification},
we prove a stronger result (as Theorem \ref{theorem-45}) which states that 
the algorithm finishes in a finite number of steps and that $W(l)$ is 
essentially $X(j)$, for some $j > l$. This will readily
imply Theorem \ref{algo-thm} as a direct consequence.

\begin{theorem}\label{algo-thm}
For general $k \in \mathbb{N}$, the sequence $X_{k}(l)$ is $2^{n_k}$-simple 
for some $n_k \in \mathbb{N}$. 
\end{theorem}

\noindent
{\bf Note}. The operators $F$ and $T$, defined in (\ref{Fdef}) and (\ref{Tdef})
respectively, play an important role in the proof of this 
conjecture.  
\\

\begin{definition}
Let $\omega(l)$ be the number of steps required for the algorithm to yield a 
constant sequence. The sequence of integers 
\begin{eqnarray}
\Omega(l) & := & \left\{ n_{1}, \, n_{2}, \, n_{3}, \cdots, n_{\omega(l)} 
\right\}
\end{eqnarray}
\noindent
is called the {\em reduction sequence} of $l$.  The number $\omega(l)$ will
be called the {\em reduction length} of $l$. The constant sequence obtained 
after $\omega(l)$ steps is called the {\em reduced constant}. 
\end{definition}

In Corollary \ref{omegal} we enumerate $\omega(l)$ as the number of ones in 
the binary expansion of $l$. Therefore the algorithm yields a 
constant sequence in a finite number of steps.  In fact, the 
algorithm terminates in $O(\log_{2}(l))$ steps. \\

Table 1 shows the results of the algorithm for $4 \leq l \leq 15$.

\begin{table}
\caption{Reduction sequence for $1 \leq l \leq 15$.}
\begin{tabular}{||c|c|c||}
\hline
$l$ & \text{binary form } & $\Omega(l) $ \\
\hline
4  & 100 & 3 \\
5 & 101 & 1, \, 2 \\
6 & 110 & 2, \, 1 \\
7 & 111 & 1, \, 1, \, 1 \\
8 & 1000 & 4 \\
9 & 1001 & 1, \, 3 \\
10 & 1010 & 2, \, 2 \\
11 & 1011 & 1, \, 1, \, 2 \\
12 & 1100 & 3, \, 1 \\
13 & 1101 & 1, \, 2, \, 1 \\
14 & 1110 & 2, \,1, \,1 \\
15 & 1111 & 1, \,1, \,1, \,1 \\
\hline
\end{tabular}
\end{table}

\medskip

We also provide a combinatorial interpretation of $\Omega(l)$. This requires 
the composition of the index $l$. 

\begin{definition}
\label{def-compo}
Let $l \in \mathbb{N}$. The {\em composition} of $l$, denoted 
by $\Omega_{1}(l)$, is 
defined as follows: write $l$ in binary form. 
Read the sequence from right to left. The first part of $\Omega_{1}(l)$ is 
the number of digits up to and including the first $1$ read in the 
corresponding binary sequence; the second one is the number of additional 
digits up to and including the second $1$ read, and so on. 
\end{definition}

For example, 
\begin{equation}
\Omega_{1}(13) = \{ 1, \, 2, \, 1 \} \text{ and } 
\Omega_{1}(14) = \{ 2, \, 1, \, 1 \}. 
\end{equation}

\medskip

Observing the values in Table 1, $\Omega_{1}(13) = \Omega(13)$ and 
$\Omega_{1}(14) = 
\Omega(14)$. We claim that this is always true.  \\

\medskip
\begin{theorem}
\label{mainthm}
The reduction sequence $\Omega(l)$ associated to an integer $l$ 
is the sequence
of compositions of $l$, that is, 
\begin{equation}
\Omega(l) = \Omega_{1}(l)
\end{equation}
\end{theorem}

\noindent

This assertion is slighted restated and proved in Section \ref{sec-algo},
as \ref{thm-reduc}. See the Note following the latter theorem.

\section{The $2$-adic valuations of $A_{l,m}$} \label{sec-valuation} 
\setcounter{equation}{0}

We now present two proofs of theorem \ref{2adicall}. 

\medskip

\begin{proof}
\noindent
{\bf First proof}. Define the numbers
\begin{eqnarray}
B_{l,m} & := & \frac{A_{l,m}}{2^{l} (m+1-l)_{2l}}.
\end{eqnarray}
\noindent
We need to prove that $B_{l,m}$ is odd. The WZ-method \cite{aequalsb}
shows that the numbers $B_{l,m}$ satisfy the
recurrence
\begin{eqnarray}
B_{l-1,m} & = & (2m+1)B_{l,m} -(m-l)(m+l+1)B_{l+1,m}, \quad 1 \leq l \leq 
m-1. \nonumber
\end{eqnarray}
\noindent
The initial values $B_{m,m} = 1$ and $B_{m-1,m} = 2m+1$ show that $B_{l,m}$ 
is an odd integer as required. \\

\noindent
{\bf Second proof}. We have 
\begin{eqnarray}
\nu_{2} \left( A_{l,m} \right) & = & l + 
\nu_{2} \left( \sum_{k=l}^{m} T_{m,k} 
\, \frac{(m+k)!}{(m-k)! \, (k-l)!} \right),
\label{identity1}
\end{eqnarray}
\noindent
where 
\begin{equation}
T_{m,k} = \frac{(2m-2k)!}{{2^{m-k} \, (m-k)!}}.
\end{equation}
\noindent
The identity
\begin{equation}
T_{m,k} = \frac{(2(m-k))!}{2^{m-k} \, (m-k)!} = 
(2m-2k-1)(2m-2k-3) \cdots 3 \cdot 1 
\end{equation}
\noindent
shows that $T_{m,k}$ is an odd integer. Then
(\ref{identity1}) can be written as
\begin {eqnarray}
\nu_{2}(A_{l,m}) & = & l + \nu_{2} \left( \sum_{k=0}^{m-l} 
T_{m,l+k} \, \frac{(m+k+l)!}{(m-k-l)! \, k!} \right) \nonumber \\
 & = & l + \nu_{2} \left( \sum_{k=0}^{m-l} T_{m,l+k} 
\frac{(m-k-l+1)_{2k+2l} }{k!} \right). \nonumber 
\end{eqnarray}
The term corresponding to $k=0$ is singled out as we write
\begin {eqnarray}
\nu_{2}(A_{l,m}) & = & l + \nu_{2} \left( T_{m,l} (m-l+1)_{2l} + 
\sum_{k=1}^{m-l} T_{m,l+k} 
\frac{(m-k-l+1)_{2k+2l} }{k!} \right). \nonumber 
\end{eqnarray}
\noindent
The claim
\begin{eqnarray}
\nu_{2} \left( \frac{(m-k-l+1)_{2k+2l}}{k!} \right) & > & 
\nu_{2} ( (m-l+1)_{2l}) \label{claim}
\end{eqnarray}
\noindent
for any $k, \, 1 \leq k \leq m-l$, will complete the proof.

To prove (\ref{claim}) we use the identity
\begin{eqnarray}
\frac{(m-k-l+1)_{2k+2l}}{k!} & = & 
(m-l+1)_{2l} \cdot \frac{(m-l-k+1)_{k} \, (m+l+1)_{k}}{k!} \nonumber
\end{eqnarray}
\noindent
and the fact that the product of $k$ consecutive 
numbers is always divisible by $k!$. This follows from the identity
\begin{equation}
\frac{(a)_{k}}{k!} = \binom{a+k-1}{k}. 
\end{equation}
\noindent
Now if $m+l$ is odd,
\begin{equation}
\nu_{2} \left( \frac{(m-l-k+1)_{k}}{k!} \right) \geq 0 \text{  and }
\nu_{2}((m+l+1)_{k}) > 0, 
\end{equation}
\noindent
and if $m+l$ is even
\begin{equation}
\nu_{2} \left( \frac{(m+l+1)_{k}}{k!} \right) \geq 0 \text{  and }
\nu_{2}((m-l-k+1)_{k}) > 0. 
\end{equation}
\noindent
This proves (\ref{claim}) and establishes the theorem.
\end{proof}

\section{Properties of the function $\nu_{2}(A_{l,m})$} \label{sec-graph} 
\setcounter{equation}{0}

In this section, we describe properties of the function 
$\nu_{2}(A_{l,m})$ for $l$ fixed and $m \geq l$. In particular, we show that 
each of these sequences has a block structure. 

\begin{theorem}
\label{thm-jump}
Let $l \in \mathbb{N}$ be fixed. Then for $m \geq l$, we have
\begin{eqnarray}
\nu_{2}(A_{l,m+1}) - \nu_{2}(A_{l,m}) & = & \nu_{2}(m+l+1) - \nu_{2}(m-l+1). 
\label{jump}
\end{eqnarray}
\end{theorem}
\begin{proof}
From (\ref{2valuel}) and $(a)_{k} = (a+k-1)!/(a-1)!$, we have
\begin{equation}
\nu_{2}(A_{l,m}) = \nu_{2} \left( \frac{(m+l)!}{(m-l)!} \right) + l. 
\end{equation}
\noindent
This implies 
\begin{eqnarray}
\nu_{2}(A_{l,m+1}) - \nu_{2}(A_{l,m}) & = & 
\nu_{2} \left( \frac{(m+l+1)!}{(m-l+1)!} \right) - 
\nu_{2} \left( \frac{(m+l)!}{(m-l)!} \right)  \nonumber \\
 & = & \nu_{2} \left( \frac{(m+l+1)! \, (m-l)!}{(m-l+1)! \, (m+l)!} 
\right) \nonumber \\
 & = & \nu_{2} \left( \frac{m+l+1}{m-l+1} \right). \nonumber
\end{eqnarray}
\noindent
The result follows from here. 
\end{proof}

The next corollary is a  special case of Theorem \ref{thm-jump}. \\

\begin{corollary}
The sequence $\nu_{2}(A_{l,m})$ satisfies 

\noindent
$1) \, \nu_{2}(A_{l,l+1}) = \nu_{2}(A_{l,l}). $

\noindent
$2)$ For $l$ even, 
\begin{equation}
\nu_{2}(A_{l,l+3}) = \nu_{2}(A_{l,l+2}) = 
\nu_{2}(A_{l,l+1}) = \nu_{2}(A_{l,l}). 
\nonumber
\end{equation}

\noindent
$3)$ The sequence $\nu_{2}(A_{1,m})$ is $2$-simple, i.e., 
$\nu_{2}(A_{1,m+1}) = \nu_{2}(A_{1,m})$. In fact, 
\begin{equation}
A_{1,m} = \{ 2, \, 2, \, 3, \, 3, \, 2, \, 2, \, 4, \, 4, \, 2, \, 2, \ldots \}.
\nonumber
\end{equation}
\end{corollary}

\medskip

Fix $k, \, l \in \mathbb{N}$ and let $\mu: = 1 + \nu_{2}(l)$. Define the 
sets
\begin{equation}
C_{k,l} := \{ l + k \cdot 2^{\mu} +j: \, 0 \leq j \leq 2^{\mu} -1  \, \}. 
\end{equation}
\noindent
Clearly the cardinality of $C_{k,l}$ is $2^{\mu}$. For 
example, if $l \in \mathbb{N}$ 
is odd, then $\mu=1$ and 
\begin{equation}
C_{k,l} = \{ l+2k, \, l+2k+1 \}. 
\end{equation}

The next result is immediate.

\begin{lemma}
\label{lemma1}
The sets $\{ C_{k,l}: \, k \geq 0\}$ form  a disjoint 
partition; namely,
\begin{equation}
\{ m \in \mathbb{N}: \, m \geq l \} = \bigcup_{k \geq 0} C_{k,l}, 
\end{equation}
\noindent
and $C_{r,l} \cap C_{t,l} = \emptyset$, whenever $r \neq t$. 
\end{lemma}

\begin{lemma}
\label{lemma2}
Fix $l \in \mathbb{N}$.

\noindent
$1)$ The sequence $\{ \nu_{2}(A_{l,m}): \, 
m \in C_{k,l} \, \}$ is constant. 
We denote this value by $\nu_{2}(C_{k,l})$. 

\noindent
$2)$ For $k \geq 0$, $\nu_{2}(C_{k+1,l}) \neq \nu_{2}(C_{k,l})$. 
\end{lemma}
\begin{proof}
Suppose $0 \leq j \leq 2^{\mu}-2$. Then 
\begin{equation}
\nu_{2}(2l+k \cdot 2^{\mu}) 
\geq \nu_{2}(k \cdot 2^{\mu}) \geq \mu > \nu_{2}(j+1), 
\end{equation}
\noindent
and hence 
\begin{equation}
\nu_{2}(2l+k \cdot 2^{\mu} + j+1) = \nu_{2}(j+1) = 
\nu_{2}(k \cdot 2^{\mu} + j+1).
\end{equation}
\noindent
Using these facts and (\ref{jump}), we obtain
\begin{eqnarray}
\nu_{2}(A_{l,l+k \cdot 2^{\mu}+j+1}) - \nu_{2}(A_{l,l+k \cdot 2^{\mu} +j}) 
 & = & \nu_{2}(2l+k \cdot 2^{\mu} +j+1) - 
\nu_{2}(k \cdot 2^{\mu} + j + 1) \nonumber \\
& = & \nu_{2}(j+1) - \nu_{2}(j+1) = 0 \nonumber 
\end{eqnarray}
\noindent
for consecutive values in  $C_{k,l}$. This proves 
part $1)$. To prove part $2),$ it 
suffices to take elements $l + k \cdot 2^{\mu} + 2^{\mu}-1 \in C_{k,l}$ and 
 $l + (k+1) \cdot 2^{\mu} \in C_{k+1,l}$ and compare their $2$-adic 
values. Again by 
(\ref{jump}), we have
\begin{eqnarray}
\nu_{2}(A_{l,l+(k+1) \cdot 2^{\mu}}) - 
\nu_{2}(A_{l,l+(k+1) \cdot 2^{\mu}-1}) & = & 
\nu_{2}(2l+(k+1) \cdot 2^{\mu}) - \nu_{2}((k+1) \cdot 2^{\mu}) \nonumber  \\
& = & \mu + \nu_{2}(2l \cdot 2^{-\mu} + k+1) - \mu - \nu_{2}(k+1) \nonumber \\
& = & \nu_{2}(2l \cdot 2^{-\mu} + k+1) - \nu_{2}(k+1) \neq 0. \nonumber 
\end{eqnarray}
\noindent
The last step follows from $2l \cdot 2^{-\mu}$ being odd and  thus 
$2l \cdot 2^{-\mu} + k+1$ and $k+1$ having opposite parities. This 
completes the proof. 
\end{proof}

\medskip

\begin{theorem}
\label{period-thm}
For each $l \geq 1$, the set $\{ \nu_{2}(A_{l,m}): \, m \geq l \, \}$ is an
$s$-simple sequence, with $s = 2^{1+ \nu_{2}(l)}$. 
\end{theorem}
\begin{proof}
From Lemma \ref{lemma1} and Lemma \ref{lemma2}, we know that $\nu_{2}( \cdot)$ 
maintains 
a constant value on each of the disjoint sets $C_{k,l}$. The length of 
each of these blocks is $2^{1 + \nu_{2}(l)}$. 
\end{proof}

\medskip

\section{The algorithm and its combinatorial interpretation}
\label{sec-algo} 
\setcounter{equation}{0}

\noindent
The proof of the Theorem \ref{mainthm} requires some preliminaries. \\

\noindent
A) Given the values of $\Omega_{1}(l)$ for 
$2^{j} \leq l \leq 2^{j+1}-1$, the list for
$2^{j+1} \leq l \leq 2^{j+2}-1$ is formed according to the following rule: \\

\noindent
$l$ is even: add $1$ to the first part of $\Omega_{1}(l/2)$ to obtain 
$\Omega_{1}(l)$; \\

\noindent
$l$ is odd: prepend a $1$ to $\Omega_{1} \left( \frac{l-1}{2} \right)$ to obtain 
$\Omega_{1}(l)$. \\

\noindent
This is clear: if $x_{1}x_{2} \cdots x_{t}$ is the binary representation 
of $l$, then $x_{1}x_{2} \cdots x_{t}0$ is the one  for $2l$. Thus, the first 
part of $\Omega_{1}(2l)$ is 
increased by $1$, due to the extra $0$ on the right. The
relative position of the remaining $1s$ stays the same. A similar argument 
takes care of $\Omega_{1}(2l+1)$. The extra $1$ that is placed at the end of the
binary representation gives the first $1$ in $\Omega_{1}(2l+1)$. \\

\noindent
B) We now relate the $2$-adic valuation of $A_{l,m}$ to 
that of $A_{\lfloor{l/2 \rfloor},m}$. \\

\begin{proposition}
\label{2adicalm1}
Let 
\begin{equation}
\lambda_{l}: = \frac{1- (-1)^{l}}{2}, \quad M_{0}: = \lfloor{ \frac{m + 
\lambda_{l}}{2}  \rfloor}.
\label{evenodd}
\end{equation}
\noindent
Then
\begin{equation}
\nu_{2}(A_{l,m}) = 2l - \lfloor{ l/2 \rfloor} + 
\lambda_{l} \nu_{2} ( M_{0} - \lfloor{ l/2 \rfloor} ) + 
\nu_{2}( A_{\lfloor{l/2 \rfloor}, M_{0}}). 
\end{equation}
\end{proposition}
\begin{proof}
We present the details for $\nu_{2}(A_{2l,2m})$.
Theorem \ref{2adicall} gives
\begin{eqnarray}
\nu_{2}(A_{2l,2m}) &  = &  \nu_{2}( (2m-2l+1)_{4l}) + 2l \nonumber \\
 & = & \nu_{2}( (2m-2l+1)(2m-2l+2) \cdots (2m+2l-1)(2m+2l) ) + 2l \nonumber \\
 & = & \nu_{2}( 2^{2l} (m-l+1)(m-l+2) \cdots (m+l)) + 2l \nonumber \\
 & = & 4 l + \nu_{2}( (m-l+1)_{2l}) \nonumber \\
 & = & 3l + \nu_{2}(A_{l,m}). \nonumber
\end{eqnarray}
\noindent
A repeated application deals with the general case. 
\end{proof}

\medskip

\begin{corollary}
\label{valAlm}
The $2$-adic valuation of $A_{l,m}$ satisfies
\begin{equation}
\nu_{2}(A_{l,m}) = 2l + \nu_{2}(l!) + 
\sum_{k \geq 0} \lambda_{\lfloor{ l/2^{k} \rfloor}} \, 
\nu_{2} ( M_{k} - \lfloor{l/2^{k+1} \rfloor} )
\end{equation}
\noindent
where
\begin{equation}
\label{Mkdef}
M_k=
\lfloor{
\frac{m+\lambda_l+2\lambda_{\lfloor{l/2}\rfloor}+\cdots 
+2^k\lambda_{\lfloor{l/2^k}\rfloor}}{2^{1+k}}}\rfloor=
\lfloor{
\frac{m+\sum_{n=0}^k2^n\lambda_{\lfloor{l/2^n}\rfloor}}{2^{1+k}}\rfloor}.
\end{equation}
\end{corollary}
\begin{proof}
This is  a repeated application of Proposition \ref{2adicalm1}. The first term 
results from 
\begin{eqnarray}
\sum_{k \geq 0} \left( 2 \lfloor{ \frac{l}{2^{k}} \rfloor} - 
\lfloor{ \frac{l}{2^{k+1}} \rfloor} \right) & = & 2l + \sum_{k \geq 1} 
\lfloor{ \frac{l}{2^{k}} \rfloor} \nonumber \\
 & = & 2l + \nu_{2}(l!). \nonumber  
\end{eqnarray}
\end{proof}

\medskip

\section{Verification of the Algorithm and the Reduction sequence} 
\label{verification} 
\setcounter{equation}{0}

In this section we establish Theorems \ref{algo-thm} and 
\ref{mainthm}, strengthend as 
Theorem \ref{theorem-45} and restated as Theorem \ref{thm-reduc}, 
respectively. First we 
prove that the reduction process alluded to in the Introduction 
is in fact an 
algorithm. This will be followed by a proof that the reduction sequence 
that comes from completing the algorithm on $X(l)$ is identical to the
composition sequence of the integer $l$. 

We now remind the reader of some definitions and nomenclature: $\Omega(l)$ is 
the reduction sequence of $X(l)$, and $\Omega_1(l)$ is the composition of 
the integer $l$. Also, $X_k(l)$ is the new sequence outputted at the end 
of the $(k-1)$th 
cycle of the algorithm, and we also use the previously defined operators
$$F(\{a_1,a_2,a_3,\dots\})=\{a_1,a_1,a_2,a_3,\dots\}$$
and
$$T(\{a_1,a_2,a_3,\dots\})=\{a_1,a_3,a_5,a_7,\dots\}.$$
Observe that $T(\{a_m: m\geq 1\})=\{a_{2m-1}: m\geq 1\}.$
Recall the constant sequence
$$c:=\{\nu_2(m): m\geq 1\}=\{0,1,0,2,0,1,0,3,0,\dots\}.$$ \\

\noindent
\bf Convention: \rm We write $A_{l,m}$ and $A_l(m)$ interchangeably.

\medskip

\noindent
{\bf Notations}: Bold-face letters will denote constant sequences, as in, 
\rm \bf t \rm =$\{t,t,t,\dots\}$. The initial sequence is 
$X(l)=\{\nu_2(A_l(m-1+l)): m\geq 1\}$. Note from Theorem 
\ref{2adicall} that
$$X(l)=\{\nu_2\left( 
\frac{(m-1+2l)!}{(m-1)!}\right)+l: m\geq 1\}.$$

\noindent
\begin{definition}
A sequence \bf a \rm =$\{a_1,a_2,a_3,\dots\}$ is 
a \it translate \rm of \bf b \rm =$\{b_1,b_2,b_3,\dots\}$ if 
\bf a = b + t\rm, for some constant sequence \bf t\rm .
\end{definition}

\medskip

Now, before proving the next main result, we consider the base case $l=1$.  \\

\begin{lemma}
\label{lemmaone}
The initial case $l=1$ satisfies
\begin{equation}
W(1)=F(T(X(1))-c)=\bf{2}.
\label{evenorodd} 
\end{equation}
\end{lemma}
\begin{proof}
Since $\nu_2(A_1(m))=\nu_2(m(m+1))+1$ and $\nu_2(2m-1)=0$, we have
$$T(X(1))=\{\nu_2((2m-1)(2m))+1: m\geq 1\}=
\{\nu_2(m)+2: m\geq 1\}=c+\text{\bf 2\rm}.$$
Then the assertion follows from $F( {\bf t}) = {\bf t}$ for a constant 
${\bf t}$. 
\end{proof}

Remember now that $X(l)$ is $2^n$-simple, hence so are its' translates. 
Thus, the next result will suffice to prove Theorem \ref{algo-thm}.
 
\bigskip

\begin{theorem}
\label{theorem-45}
The algorithm terminates after finite iterations. Further, in 
each cycle, $W(l)$ is a translate of $X(j)$, for some $j < l$.  
\end{theorem}
\begin{proof}
Start by rewriting the terms in $X(l)$ as
$$\nu_2\left(\frac{(m-1+2l)!}{(m-1)!}\right)+l=
\nu_2((m-1+2l)(m-2+2l)\cdots (m+1)m),\qquad m\geq 1.$$
Then, the operator $T$ acts on these to yield (for $m\geq 1$)
\begin{equation}
\nu_2((2m-2+2l)(2m-3+2l)\cdots (2m)(2m-1))+l  
\nonumber
\end{equation}
\begin{eqnarray}
& = & \nu_2((m-1+l)\cdots (m))+2l \nonumber \\
& = & \nu_2\left(\frac{(m-1+l)!}{(m-1)!}\right)+2l. 
\label{nu2long}
\end{eqnarray}

\noindent
\bf Case I: \rm $l$ is even. From (\ref{nu2long}), we 
can easily obtain the relation (with $l_2=l/2$) 
$$T(X(l))=\{\nu_2\left(\frac{(m-1+2l_2)!}{(m-1)!}\right)+l_2+t:m\geq 1\}=
X(l_2)+\text{\bf t\rm}, \qquad t=3l_2.$$

\noindent
\bf Case II: \rm $l$ is odd. Upon 
subtracting the sequence 
$c=\{\nu_2(m): m\geq 1\}$ from (\ref{nu2long}) 
and letting $l_1=(l-1)/2$, we get that 
$$\nu_2\left(\frac{(m+l-1)!}{m!}\right)+2l=
\nu_2\left(\frac{(m+2l_1)!}{m}\right)+2l=
\nu_2\left(\frac{(m+2l_1)!}{m!}\right)+l_1+(3l_1+2), 
$$
\noindent
for $m\geq 1.$ 
Finally, apply the operator $F$ to the last sequence and find
$$W(l)=\{\nu_2\left(\frac{(m-1+2l_1)!}{(m-1)!}\right)+l_1+t:m\geq 1\}=
X(l_1)+\text{\bf t\rm}, \qquad t=3l_1+2.$$
Here, we have  utilized the fact that $\nu_2(2l+(l-1)!)=\nu_2(2l+l!)=1$ 
which 
is valid for $l$ odd. This justifies that the first term augmented in the 
sequence, as a result of the action of $F$, coincides 
with the next term (these are 
values at $m=1$ and $m=2$, respectively, in $X(l_1)$). \\

\noindent
We can now conclude that in either of the two cases (or a combination 
thereof), the index $l$ shrinks dyadically as $l_{1}$ or $l_{2}$. Thus 
the reduction algorithm must end in a finite step into a 
translate of $X(1)$. Since Lemma \ref{lemmaone} 
handles $X(1)$, the proof is 
completed.
\end{proof}

\bigskip

\begin{theorem}
\label{thm-reduc}
Let $\{ k_{1}, \cdots, k_{n}: \, 0 \leq k_{1} < k_{2} < \cdots < k_{n} \}$, be
the unique collection of distinct positive integers such that
\begin{equation}
 l = \sum_{i=1}^{n} 2^{k_{i}}. 
\end{equation}
\noindent
Then the reduction sequence of $l$ is $\{ k_{1}+1, \, k_{2}-k_{1}, \cdots, 
k_{n}-k_{n-1} \}$. 
\end{theorem}

\noindent
{\bf Note}. The argument of the proof is to check that the rules of 
formation for 
$\Omega_{1}(l)$ also hold for the reduction sequence $\Omega(l)$. This 
will incidentally elaborates the connection with 
\ref{mainthm}. The proof is divided according to the parity of $l$. 

\begin{proof}
The case $l$ odd starts with $l=1$, where the block length is $2$. From 
Theorem \ref{2adicall}
we obtain a constant sequence after iterating the algorithm once. Thus the algorithm terminates and the
reduction sequence for $l=1$ is $\Omega(1) = \{ 1 \}$. \\

Now consider the general even case:  $X(2l)$.  
Applying $T$ to this sequence yields a 
translate of $X(l)$ by Theorem
\ref{theorem-45}; this does not affect the reduction 
sequence $\Omega(l)$, but the doubling 
of block length increases the first term of $\Omega(l)$ by $1$. Therefore 
\begin{equation}
\Omega(2l) = \{ k_{1}+2, \, k_{2}-k_{1}, \, \cdots, k_{n}-k_{n-1} \}. 
\end{equation}
\noindent
This is precisely what happens to the binary digits of $l$: if 
$$l = \sum_{i=1}^{n} 2^{k_{i}}, \text{ then }
2l = \sum_{i=1}^{n} 2^{k_{i}+1}. $$
\noindent
This concludes the argument for even indices. 

For the general odd case, $X(2l+1)$, we apply $T$, subtract $c$ and 
then apply $F$.  Again, by Theorem \ref{theorem-45}, this gives us a 
translate 
of $X(l)$. We conclude that, if the reduction sequence of $l$ is 
\begin{equation}
\{ k_{1}+1, \, k_{2}-k_{1}, \cdots, k_{n}-k_{n-1} \},
\end{equation}
\noindent
then that of $2l+1$ is 
\begin{equation}
\{ 1, \, k_{1}+1, \, k_{2}-k_{1}, \cdots, k_{n}-k_{n-1} \}.
\end{equation}
\noindent
This is  precisely the behavior of $\Omega_{1}$.  The proof is complete.
\end{proof}

\medskip

\begin{corollary}
\label{redvalue}
The reduced constant is $2l + \nu_{2}(l!) = \nu_{2}(A_{l,l})$. 
\end{corollary}
\begin{proof}
In Corollary \ref{valAlm}, subtract the last term as per the reduction
algorithm (or as implied by Theorem \ref{thm-reduc} or Theorem \ref{mainthm}).
\end{proof}

\begin{corollary}\label{omegal}
The set $\Omega(l)$ has cardinality
\begin{equation}
s_{2}(l) = \text{the number of ones in the binary expansion of } l. 
\end{equation}
\end{corollary}

\medskip

\noindent
{\bf Remarks}: \\

Write $l$ in
the binary form: $l=\sum_{j=1}^n2^{k_j}$ with $0\leq k_1<\cdots <k_n$. Then, 
for the $M_{k}$ defined in (\ref{Mkdef}) can be rewritten as 
$$M_{k_i}=\lfloor{\frac{m+\sum_{j=1}^i2^{k_j}}{2^{1+k_i}}}\rfloor.$$

\noindent
1) In light of this, Corollary \ref{valAlm} may be given in the form
\begin{equation}
\nu_2(A_{l,m})=2l+\nu_{2}(l!)+\sum_{i\geq 1}
\nu_2\left(M_{k_i}-\lfloor{l/2^{1+k_i}}\rfloor\right). \label{star}
\end{equation}

\medskip 

\noindent
2) Observe also that $\nu_2(M_{k_i}-\lfloor{l/2^{1+k_i}}\rfloor)$ is 
a \it $2^{1+k_i}$-simple \rm 
sequence, i.e. it has constant blocks of length $2^{1+k_i}$.

\medskip 

\noindent
3) The sequence
$\nu_2(A_{l,m})$ inherits its \it $2^{1+k_1}$-simple \rm structure from 
the term 
$\nu_2(M_{k_1}-\lfloor{l/2^{1+k_1}}\rfloor)$, which has the lowest 
period (or highest frequency) in the decomposition (\ref{star}). Notice that 
this is consistent with Theorem \ref{period-thm}, since $k_{1} = \nu_{2}(l)$.

\medskip 

\noindent
4) The sequence $(\dots, \lambda_{\lfloor{l/2}\rfloor},\lambda_{l})$ is the 
binary code for $l$, and $(\dots,k_2+1,k_1+1)$ are the exponents of $2$ in 
the binary format of $2l$. 

\medskip 

\noindent
5) For fixed $l$, we can construct the sequence $\nu_2(A_{l,m})$ by 
reversing the algorithm. Write the binary code for 
$2l=\sum_{j=1}^n2^{1+k_j}$, and 
then, starting with the $\infty$-simple (constant) sequence $3l-s_2(l)$, then 
add the $2^{1+k_1}-, 2^{1+k_2}-,\dots, 2^{1+k_n}-$simple 
sequences $\nu_2(M_{k_i}-\lfloor{l/2^{1+k_i}}\rfloor).$ Here, the successive 
differences $(1+k_j)-(1+k_{j-1})=k_j-k_{j-1}$, for $j=1=1,\dots,n$, encode 
the \it period switching-gaps \rm (or \it indices of sequence shifting \rm as 
compared to the preceding stages) on the one hand, and the \it integer 
composition \rm of $2l$ on the other. This shades more light into the 
bijective 
relationship between $\Omega(l)$ and $\Omega_1(l)$ that has been proven 
in Theorem \ref{mainthm}. 

\bigskip

\noindent
{\bf Note}. The function $s_{2}(l)$ has recently appeared in a different 
divisibility problem. In these papers it is denoted by $d(l)$. 
Lengyel \cite{lengyel1} conjectured, and De 
Wannemacker \cite{wannemacker1} proved, that the $2$-adic valuation of
the Stirling numbers of the second kind $S(n,k)$ is given by 
\begin{equation}
\nu_{2} ( S(2^{n},k)) = s_{2}(k)-1. 
\label{stirlval1}
\end{equation}
\noindent
The Stirling numbers are given by the identity
\begin{equation}
x^{n} = \sum_{k=0}^{n} S(n,k)  x(x-1)(x-2) \cdots (x-k+1)
\end{equation}
\noindent
and they count the number of ways to partition a set with $n$ elements into
exactly $k$ nonempty subsets.
De Wannemacker \cite{wannemacker2} also established the inequality 
\begin{equation}
\nu_{2}( S(n,k) ) \geq s_{2}(k) - s_{2}(n), \quad 0 \leq k \leq n.
\end{equation}
The study of the $2$-adic valuation of Stirling numbers suggests that 
\begin{equation}
\nu_{2} ( S(2^{n}+1,k+1)) = s_{2}(k) -1, 
\label{stirlval2}
\end{equation}
\noindent
which is a companion of (\ref{stirlval1}).  \\

\section{A connection with the Collatz problem}
\setcounter{equation}{0}

The numbers 
\begin{equation}
a_{m}:= \nu_{2}(A_{1,m}) -1 =  \nu_{2}(m(m+1)),
\end{equation}
\noindent
given in (\ref{2value1}), also appear in the well-known
{\em Collatz} or $3x+1$ problem. 
Define a sequence by $x_{0}(m) = m$ and let 
$x_{k+1}(m) = T(x_{k}(m))$, where
\begin{equation}
T(i)= \begin{cases}
             \frac{1}{2}i & \text{ if } i \text{ is even}, \\
             \frac{1}{2} (3i+1) & \text{ if } i \text{ is odd}.
         \end{cases}
\end{equation}
\noindent
The {\em orbit} of $m \in \mathbb{N}$ is the set 
\begin{equation}
\frak{O}(m) := \{ m, \, T(m), \, T^{2}(m), \ldots \}.
\end{equation}
\noindent
The main conjecture for this problem is that {\em every orbit ends in the
cycle} $1 \to 2 \to 1$. The reader will find 
in \cite{lagarias1} an introduction to this problem
and \cite{chamberland1,lagarias2} contain annotated bibliographies. \\

The connection with our work is given in the next theorem. \\

\begin{theorem}
Let $m \in \mathbb{N}$. Then $a_{m}:= \nu_{1}(A_{1,m})-1 = \nu_{2}(m(m+1))$ is
the first time at which the orbit $\frak{O}(m)$ changes parity. That is,
\begin{equation}
m \equiv T(m) \equiv T^{2}(m) \equiv \cdots \equiv T^{a_{m}-1}(m) \, 
\not \equiv T^{a_{m}}(m)  \,
\text{ mod } 2.
\end{equation}
\end{theorem}
\begin{proof}
Suppose $m$ is odd and write it as $m = 2^{j}n - 1$, with $n$ odd. Then 
\begin{equation}
j = \nu_{2}(m+1) \text{ and } n = \frac{m+1}{2^{j}} 
\end{equation}
\noindent 
are uniquely defined. Observe that
\begin{equation}
T(m) = T( 2^{j}n -1 ) = 3 \cdot 2^{j-1}n -1 
\nonumber
\end{equation}
\noindent
and for $i < j$, 
\begin{equation}
T^{i}(m) = T^{i}( 2^{j}n -1 ) = 3^{i} \cdot 2^{j-i}n -1.
\nonumber
\end{equation}
\noindent
Finally, 
\begin{equation}
T^{j}(m) = T^{j}(2^{j}n-1) = 3^{j}n -1. \nonumber
\end{equation}
\noindent
To complete the proof, observe that 
\begin{equation}
j = \nu_{2}(m+1) = \nu_{2}(m(m+1)) = N. 
\end{equation}
\noindent
In the case $m$ is even, write $m = 2^{t}m_{0}$, with $m_{0}$ odd. Then 
\begin{equation}
T^{i}(m) = 2^{t-i}m_{0}, \text{ for } 0 \leq i < t 
\end{equation}
\noindent
and 
\begin{equation}
T^{t}(m) = m_{0}. 
\end{equation}
\noindent
The proof is completed by noticing that 
\begin{equation}
t = \nu_{2}(m) = \nu_{2}(m(m+1)) = N. 
\end{equation}
\end{proof}

For example take $m=63$. Then
$x_{1}(63) = 95, \, x_{2}(63) = 143, \, x_{3}(63) = 215, \, x_{4}(63) = 323, 
\, x_{5}(63) = 485$, and $x_{6}(63) = 728$. Thus,
\begin{equation}
\frak{O}(63) = \{ 63, \, 95, \, 143, \, 215, \, 323, \, 485, \, 
\mathbf{728}, \ldots \}. 
\end{equation}
\noindent
It takes $6$ iterations to 
produce an even entry. Observe that $a_{63} = \nu_{2}((63)_{2}) = 6$. 

\medskip

Similarly, we have 

\begin{proposition}
The first time the orbit of $3^m-1$ changes parity is after 
$$\nu_2(3^m(3^m-1))=\nu_2(3^m-1)=\lambda_m+\nu_2(2m)=\nu_2(2^{1+\lambda_m}m)$$ 
\noindent
iterations.
\end{proposition}
\begin{proof}
Use the binomial theorem for $(2+1)^m-1$, while the generating function can 
be given by
\begin{equation}
\sum_{m\geq 1}\nu_2(3^m-1)x^m
=\frac{x^2}{1-x^2}+\sum_{k\geq 0}\frac{x^{2^k}}{1-x^{2^k}}.
\end{equation}
\end{proof}

\section{A symmetry conjecture on the graphs of $\nu_{2}(A_{l,m})$}
\setcounter{equation}{0}

The graphs of the function $\nu_{2}(A_{l,m})$, where we take every other 
$2^{1+ \nu_{2}(l)}$-element to reduce the  
repeating blocks to a single value, are shown in the next figures. We 
conjecture that these graphs have a symmetry property generated by
what we call 
an {\em initial segment}: from which the rest is
determined by adding a {\em central piece} followed by a {\em folding 
rule}. For example, in the case $l=1$, the first few
values of the reduced table  are 
$$ \{ 2, \, 3, \, 2, \, 4, \, 2, \, 3, \, 2, \, 5, \, 2, \, 3, \ldots \}. $$
\noindent

{{
\begin{figure}[ht]
\begin{center}
\centerline{\psfig{file=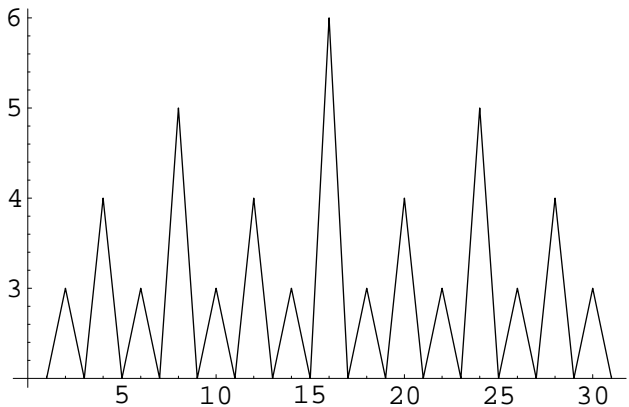,width=25em,angle=0}}
\caption{The $2$-adic valuation of $A_{1,m}$}
\label{2adic1}
\end{center}
\end{figure}
}}

The ingredients are:  \\

\noindent
{\em initial segment}: $ \{ 2, \, 3, \, 2 \} $,  \\

\noindent
{\em central piece}: the value at the 
center of the initial segment, namely $3$.  \\

\noindent
{\em rules of formation}: 
start with the initial segment and add $1$ to the 
central piece and reflect. \\

\noindent
This produces the sequence
$$\{ 2, 3, 2 \} \to \{ 2, 3, 2, 4 \} \to 
\{2,3,2,4,2,3,2\} \to 
\{ 2,3,2,4,2,3,2,5 \} \to $$
$$ \to \{2,3,2,4,2,3,2,5,2,3,2,4,2,3,2 \}. $$

\bigskip

The details are shown in Figure \ref{2adic1}.

\medskip

We have found no way to predict the initial segment
nor the central piece. Figure \ref{2adic9a} shows the
beginning of the case $l=9$. From here one could be tempted to 
anticipate that this graph extends as in the case $l=1$. This is not 
correct however, as can be seen in Figure \ref{2adic9b}. 
In fact, the initial segment is depicted
in Figure \ref{2adic9b} and its extension is shown in Figure \ref{2adic9c}.

{{
\begin{figure}[ht]
\begin{center}
\centerline{\psfig{file=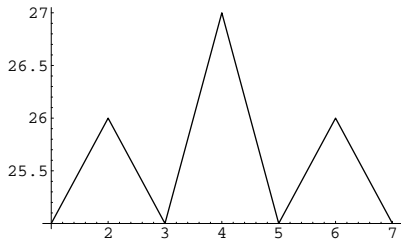,width=15em,angle=0}}
\caption{The beginning for $l=9$}
\label{2adic9a}
\end{center}
\end{figure}
}}

\medskip

{{
\begin{figure}[ht]
\begin{center}
\centerline{\psfig{file=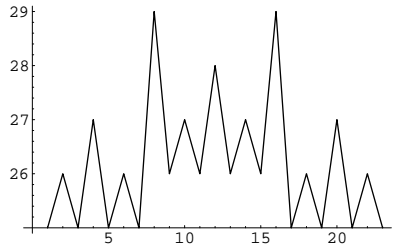,width=15em,angle=0}}
\caption{The continuation  of $l=9$}
\label{2adic9b}
\end{center}
\end{figure}
}}

\medskip

{{
\begin{figure}[ht]
\begin{center}
\centerline{\psfig{file=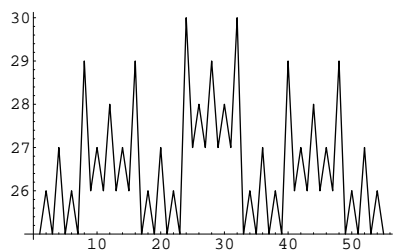,width=15em,angle=0}}
\caption{The pattern for $l=9$ persists}
\label{2adic9c}
\end{center}
\end{figure}
}}

\medskip

The initial pattern can be quite elaborate. Figure \ref{2adic53} illustrates
the case $l=53$.

{{
\begin{figure}[ht]
\begin{center}
\centerline{\psfig{file=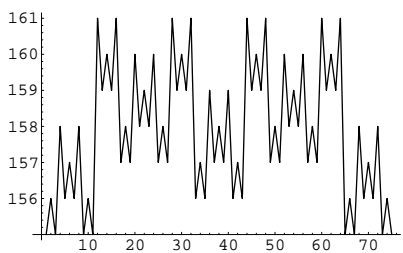,width=15em,angle=0}}
\caption{The initial pattern for $l=53$ }
\label{2adic53}
\end{center}
\end{figure}
}}

\medskip

\noindent
{\bf Acknowledgements}. The last author acknowledges the partial support of 
NSF-DMS 0409968. The second author was partially supported as  a graduate 
student by the same grant. The work of the first author was done while 
visiting Tulane University in the Spring of 2006. The authors wish to thank
Marc Chamberland for information on the $3x+1$ problem 
and Aaron Jaggard for identifying our data with the composition sequence.

\end{document}